\documentclass[final,12pt,nonatbib]{elsarticle}

\makeatletter
\let\c@author\relax
\makeatother

\usepackage{xcolor,soul,cancel}
\usepackage[color]{showkeys}

\definecolor{refkey}{rgb}{0,1,1}
\definecolor{labelkey}{rgb}{1,0,0}
\usepackage{amssymb}
\usepackage{amsmath}
\usepackage{amsthm}
\usepackage{graphicx}
\usepackage{float}
\usepackage{xcolor}
\usepackage{cases}
\usepackage{color}
\usepackage{float}
\usepackage{hyperref}
\hypersetup{
    colorlinks=true, 
    linkcolor=blue, 
    urlcolor=red, 
    linktoc=all 
}
\journal{arXiv}

\newtheorem{thm}{Theorem}
\newtheorem{lem}{Lemma}

\newtheorem{prop}[thm]{Proposition}

\newcommand{\eq} [1] {\begin{equation}\label{#1}\quad}
\newcommand{\en} {\end{equation}}
\newcommand{\scal}[1]{\langle#1\rangle}
\newcommand{\norm}[1]{\left\Vert#1\right\Vert}
\newcommand{\abs}[1]{\left\vert#1\right\vert}

\newcommand{\C}{\mathbb C}

\newcommand{\R}{\mathbb R}

\newcommand{\im}{\operatorname{Im}}

\newcommand{\re}{\operatorname{Re}}

\begin{document}

\begin{frontmatter}

\title{Numerical ranges of Foguel operators revisited. \tnoteref{support}}


\author[nyuad]{Muyan Jiang}
\ead{muyan_jiang@berkeley.com, mj2259@nyu.edu, jimmymuyan@163.com}
\author[nyuad]{Ilya M. Spitkovsky}
\address[nyuad]{Division of Science and Mathematics, New York  University Abu Dhabi (NYUAD), Saadiyat Island,
P.O. Box 129188 Abu Dhabi, United Arab Emirates}
\ead{imspitkovsky@gmail.com, ims2@nyu.edu, ilya@math.wm.edu}
\tnotetext[support]{The second author [IMS]  was supported in part by Faculty Research funding from the Division of Science and Mathematics, New York University Abu Dhabi.}

\begin{abstract}
The Foguel operator is defined as $F_T=\begin{bmatrix}S^* & T \\ 0 & S\end{bmatrix}$, where $S$ is the right shift on a Hilbert space $\mathcal H$ and $T$ can be an arbitrary bounded linear operator acting on $\mathcal H$. Obviously, the numerical range $W(F_0)$ of $F_T$ with $T=0$ is the open unit disk, and it was suggested by Gau, Wang and Wu in their LAA'2021 paper that $W(F_{aI})$ for non-zero $a\in\C$ might be an elliptical disk. In this paper, we described $W(F_{aI})$ explicitly and, as it happens, it is not.

\end{abstract}


\end{frontmatter}

\section{Introduction and the main result} 
Let $\mathcal H$ be a Hilbert space. Denoting by $\scal{.,.}$ the scalar product on $\mathcal H$ and by $\norm{.}$ the norm associated with it, recall that the {\em numerical range} of a bounded linear operator $A$ acting on $\mathcal H$ is the set 
\[ W(A):=\{ \scal{Ax,x}\colon x\in \mathcal H, \norm{x}=1\}. \]
Clearly, $W(A)$ is a bounded subset of the complex plane $\C$; the {\em numerical radius} $w(A)$ defined as 
\[ w(A)=\sup\{\abs{z}\colon z\in W(A) \}\] does not exceed $\norm{A}$ due to the Cauchy-Schwarz inequality. Furthermore, according to the celebrated Toeplitz-Hausdorff theorem, this set is also convex. The point spectrum (= the set of the eigenvalues) $\sigma_p(A)$ is contained in $W(A)$, while the whole spectrum $\sigma(A)$ is contained in its closure.  A detailed up-to-date treatment of the numerical range properties and the history of the subject can be found in a recent comprehensive monograph \cite{GauWu}. 

In this paper, we consider {\em Foguel operators}, i.e., operators acting on $\mathcal H\oplus\mathcal H$ according to the block matrix representation 
\eq{ft} F_T=\begin{bmatrix} S^* & T \\ 0 & S\end{bmatrix}.\en
Here $S$ is the right shift while $T$ is an arbitrary operator on $\mathcal H$. Since $W(S)=W(S^*)$ coincides with the open unit disk $\mathbb D$, for any $T$ we have \eq{wcd} W(F_T)\supseteq\mathbb D.\en  In particular, $w(F_T)\geq 1$. 

In \cite{GauWangWu21} some more specific estimates for $w(F_T)$ were established for an arbitrary $T$, while for $T=aI$ the exact formula $w(F_{aI})=1+\abs{a}/2$ was obtained. In the latter case it was also proved in \cite{GauWangWu21} that $W(F_{aI})$ is open, and conjectured that it is an elliptical disk with the minor half-axis $\sqrt{1+\abs{a}^2/4}$. (The case $a=0$ is of course obvious, since $F_0=S^*\oplus S$ and so $W(F_0)=\mathbb D$.) 

As we will show, this conjecture does not hold. The numerical range of $F_{aI}$ can be described explicitly, and this description is as follows.
\begin{thm}\label{th:nrd}Let $A$ be given by \eqref{ft} with $T=aI$. Then the numerical range of $A$ is symmetric with respect to both coordinate axes, bounded on the right/left by arcs of the circles centered at $(\pm 1,0)$ and having radius $r=\abs{a}/2$, and from above/below by arcs of the algebraic curve defined by the equation 
\eq{arc}\begin{gathered}
16 r^6 (u+v)^2-8 r^4 \left(u^3+(v-1)(4 u^2+5uv-u+2v^2)-v\right)\\
+r^2 \left(((u-20) u-8) v^2+2 ((u-15) u-4) (u-1) v\right. \\ \left. +((u-10) u+1) (u-1)^2\right) +(u-1)^3 (u+v-1) = 0
\end{gathered} \en
in terms of $u=x^2,\ v=y^2$. The switching points between the arcs are located on the supporting lines of $W(A)$ forming the angles $\pm\cos^{-1}\frac{\sqrt{4+r^2}-r}{2}$ with the $y$-axis.   \end{thm} 

The proof of Theorem~\ref{th:nrd} is split between the next two sections. In Section~\ref{s:proof}, the parametric description of the supporting lines to $W(F_{aI})$ is derived. Section~\ref{s:points} is devoted to the transition from that to the point equation of $\partial W(F_{aI})$. Note that the formulas in Section~\ref{s:proof} are obtained via a somewhat unexpected application of Toeplitz operators technique.

\section{$W(F_{aI})$ in terms of its supporting lines} \label{s:proof} 
As any convex subset of $\C$, the numerical range of a bounded linear operator $A$ is completely determined by the family of its supporting lines. The latter has the form 
\eq{sl} \omega\left(\lambda_{\max}(\re(\overline{\omega}A))+i\R\right), \quad \omega=e^{i\theta}\text{ with } \theta\in[-\pi,\pi].  \en 
(We are using the standard notation $\re X:=(X+X^*)/2$ for the hermitian part of an operator $X$, and $\lambda_{\max}(H)$ for the rightmost point of $\sigma(H)$ when $H$ is a hermitian operator.) 

For $A=F_T$ given by \eqref{ft} we are therefore interested in the non-invertibility of the operator 
\eq{xo} \begin{bmatrix}X_\omega-2\lambda I & \overline{\omega}T \\ 
\omega T^* & X_{\overline{\omega}}-2\lambda I\end{bmatrix}, \en where \eq{xom} X_\omega:=\omega S+\overline{\omega}S^*.\en  Due to \eqref{wcd}, we know a priori that the values of $\lambda_{\max}$ in \eqref{sl} are all greater than or equal to one. So, we may suppose that in \eqref{xo} $\lambda>1$.

Under this condition, the diagonal blocks of \eqref{xo} are invertible and, according to the Schur complement  the operator \eqref{xo} is (or is not) invertible simultaneously with $X_{\overline{\omega}}-2\lambda I - T^*(X_\omega-2\lambda I)^{-1}T$. 
For $T=aI$ this expression simplifies further to \eq{sc} X_{\overline{\omega}}-2\lambda I - \abs{a}^2(X_\omega-2\lambda I)^{-1}.\en  
Replacing $\abs{a}/2$ by $r$, as in the statement of the theorem, and multiplying \eqref{sc} by $X_\omega-2\lambda I$, we see that for $\lambda>1$ the operator \eqref{xo} is invertible only simultaneously with 
\eq{scm} X_{\overline{\omega}}X_\omega-2\lambda(X_\omega+X_{\overline{\omega}})+4(\lambda^2-r^2)I. \en 
\begin{lem} \label{th:scm}Let $\omega\neq\pm 1,\pm i$. Then the operator \eqref{scm} is not invertible if and only if $2(r^2-\lambda^2)$
lies in the range of the function 
\eq{f} f_{\lambda,\omega}(t)=\re(t^2)+\re(\omega^2)-4\lambda\re t\re\omega, \quad t\in\mathbb T. \en \end{lem}
Proof of this Lemma uses some notions concerning Hardy spaces and Toeplitz operators; interested readers are referred to excellent monographs available on the subject, e.g., \cite{BSil06} or \cite{Nik19}. 
\begin{proof}Let us consider the realization of $\mathcal H$ as the Hardy space $H^2$ of functions analytic on $\mathbb D$ with their Taylor coefficients forming an $\ell^2$ sequence. This space identifies naturally with a subspace of $L^2$ on the boundary of $\mathbb D$, the unit circle $\mathbb T$. Denote by $P$ the orthogonal projection of $L^2(\mathbb T)$ onto $H^2$, and recall that for any $\phi\in L^\infty(\mathbb T)$ the Toeplitz operator with the symbol $\phi$ is defined as
\[ (T_\phi =:)\ P\phi P\colon H^2\longrightarrow H^2. \] 
When $\phi\in H^\infty$ (the space of bounded analytic function on $\mathbb D$), $T_\phi$ is simply the multiplication by $\phi$. Such is, in particular, the shift operator $S$, which for $\mathcal H=H^2$ is nothing but the multiplication by the variable $z$:
\[ (Sf)(z)= zf(z), \quad f\in H^2. \] 
Accordingly, the operator $X_\omega$ defined by \eqref{xom} is the Toeplitz operator with the symbol $\omega t+\overline{\omega}t^{-1}$, and \eqref{scm} takes the form 
\[ P(\overline{\omega}t+\omega t^{-1})P(\omega t+\overline{\omega}t^{-1})P-2\lambda P(\omega t+\overline{\omega}t^{-1}+\overline{\omega}t+\omega t^{-1})P+4(\lambda^2-r^2)I. \] 

Replacing the middle $P$ in the first term by $I-Q$, where $Q$ is the projection of $L^2(\mathbb T)$ onto the orthogonal complement of $H^2$, we can further rewrite \eqref{scm} as  
$2T_f-\overline{\omega}^2P_0+4(\lambda^2-r^2)I$. Here $f=f_{\lambda,\omega}$ is given by \eqref{f},
and $P_0=P(t+\omega^2 t^{-1})Q(\omega^2 t+t^{-1})P$. Considering that $Pt^{-1}Q=QtP=0$, $P_0$ simplifies to $PtQt^{-1}P$, which is nothing but the rank one orthogonal projection mapping each $f\in H^2$ to its constant term. 
So, the operator \eqref{scm} is not invertible if and only if 
\eq{specond} 2(r^2-\lambda^2)\in\sigma\left(T_f-\frac{\overline{\omega}^2}{2} P_0\right). \en
Observe that the essential spectrum of the operator $T_f-\frac{\overline{\omega}^2}{2} P_0$ is the same as that of the Toeplitz operator $T_f$. Since the function $f$ is real-valued, the latter coincides with $\sigma(T_f)$. So, the right hand side of \eqref{specond} is $\sigma(T_f)$ possibly united with some isolated eigenvalues of $T_f-\frac{\overline{\omega}^2}{2} P_0$.

Our next step is to observe that for non-real $\omega^2$ such eigenvalues, if exist, cannot be real. Indeed, denoting by $\xi$ the respective eigenfunction and writing $T_f-\frac{\overline{\omega}^2}{2} P_0$ as $H+iK$ with $H,K$ hermitian, we would then have
$\scal{K\xi,\xi}=0$. But $K=\frac{\im(\omega^2)}{2}P_0$, and so $\scal{K\xi,\xi}$ is a non-zero multiple of $\abs{\xi(0)}^2$. From here, $P_0\xi=\xi(0)=0$, and $\xi$ is an eigenfunction of $T_f$ corresponding to the same eigenvalue $\mu$. This is a contradiction. 

So, for $\omega\neq \pm 1,\pm i$ condition \eqref{specond} is equivalent to $2(r^2-\lambda^2)\in\sigma(T_f)$. It remains to invoke the fact that, since the function $f$ is continuous, $\sigma(T_f)$ equals the range of $f$.
\end{proof} 
Recall that we are interested in the maximal values of $\lambda$ satisfying the conditions of Lemma~\ref{th:scm}. The explicit formulas for them are provided in the following statement.
\begin{lem}\label{th:lmax}Let $f$ be defined by \eqref{f}. Then the maximal value of $\lambda$ satisfying 
\eq{condonl} 2(r^2-\lambda^2)\in f(\mathbb T) \en 
is \eq{lmax}\lambda_{\max}(\theta)= \begin{cases} r+\abs{\re\omega} & \text{ if } 
\abs{\re\omega}\geq (\sqrt{4+r^2}-r)/2, \\
	\sqrt{1+(r/\im\omega)^2} &  {\text otherwise.} \end{cases} \en
\end{lem} 
\begin{proof}
Formulas \eqref{lmax} are invariant under substitutions $\omega\mapsto -\omega$ and $\omega\mapsto\overline{\omega}$. Since also 
$f_{\lambda,\omega}=f_{\lambda,\overline{\omega}}$ and $f_{\lambda,\omega}(t)=f_{\lambda,-\omega}(-t)$, it suffices to consider $\omega$ lying in the first quadrant only. So, in what follows $\omega=e^{i\theta}$ with $\cos\theta,\sin\theta\geq 0$. 

A straightforward calculus application yields 
\[ f(\mathbb T)=\begin{cases} [2\cos^2\theta-4\lambda\cos\theta, 2\cos^2\theta+4\lambda\cos\theta] & \text{ if } \lambda\cos\theta\geq 1, \\ [2(1-\lambda^2)\cos^2\theta-2, 2\cos^2\theta+4\lambda\cos\theta] & \text{ otherwise}. \end{cases} \]	
After some simplifications, \eqref{condonl} can therefore be rewritten as 
\[ r\in [\lambda-\cos\theta, \lambda+\cos\theta], \quad \lambda\cos\theta\geq 1, \]
or 
\[ r\in [\sqrt{\lambda^2-1}\sin\theta, \lambda+\cos\theta], \quad \lambda\cos\theta\leq 1. \]
In terms of $\lambda$, these conditions are equivalent to $\lambda\in\Lambda_1\cup\Lambda_2$, where 
\[\Lambda_1=\left[\max\{r-\cos\theta,\sec\theta\}, r+\cos\theta\right] \] and  \[ \Lambda_2=\left[r-\cos\theta,\min\{\sqrt{1+(r/\sin\theta)^2},\sec\theta\}\right], \] respectively. 

The remaining reasoning depends on the relation between $\sec\theta-\cos\theta$ and $r$ (equivalently, between $\cos\theta$ and $(\sqrt{4+r^2}-r)/2$).

\noindent 
{\sl Case 1.} $r\geq\sec\theta-\cos\theta$. Then $\max\Lambda_1=r+\cos\theta\geq\sec\theta\geq\max\Lambda_2$.

\noindent 
{\sl Case 2.} $r<\sec\theta-\cos\theta$. This inequality implies $\sqrt{1+(r/\sin\theta)^2}<\sec\theta$, and so 
$\max\Lambda_2= \sqrt{1+(r/\sin\theta)^2}$ while $\Lambda_1=\emptyset$. 

These findings agree with \eqref{lmax} thus completing the proof. \end{proof}
So, at least for $\omega\neq\pm 1,\pm i$, the supporting lines of $W({F_{aI}})$ are given by 
\eq{sup} \omega\left(\lambda_{\max}(\theta)+i\R\right), \quad \omega=e^{i\theta}\text{ with } \theta\in[-\pi,\pi]  \en 	
and $\lambda_{\max}(\theta)$ as in \eqref{lmax}. 	
Due to the continuity of $\lambda_{\max}$ at the a priori excluded values of $\omega$, this description actually holds throughout.

To illustrate, below is the plot of the supporting lines of $W(F_{aI})$ corresponding to $a=1$.

\begin{figure}[H]
    \centering
    \includegraphics[width=0.6\textwidth]{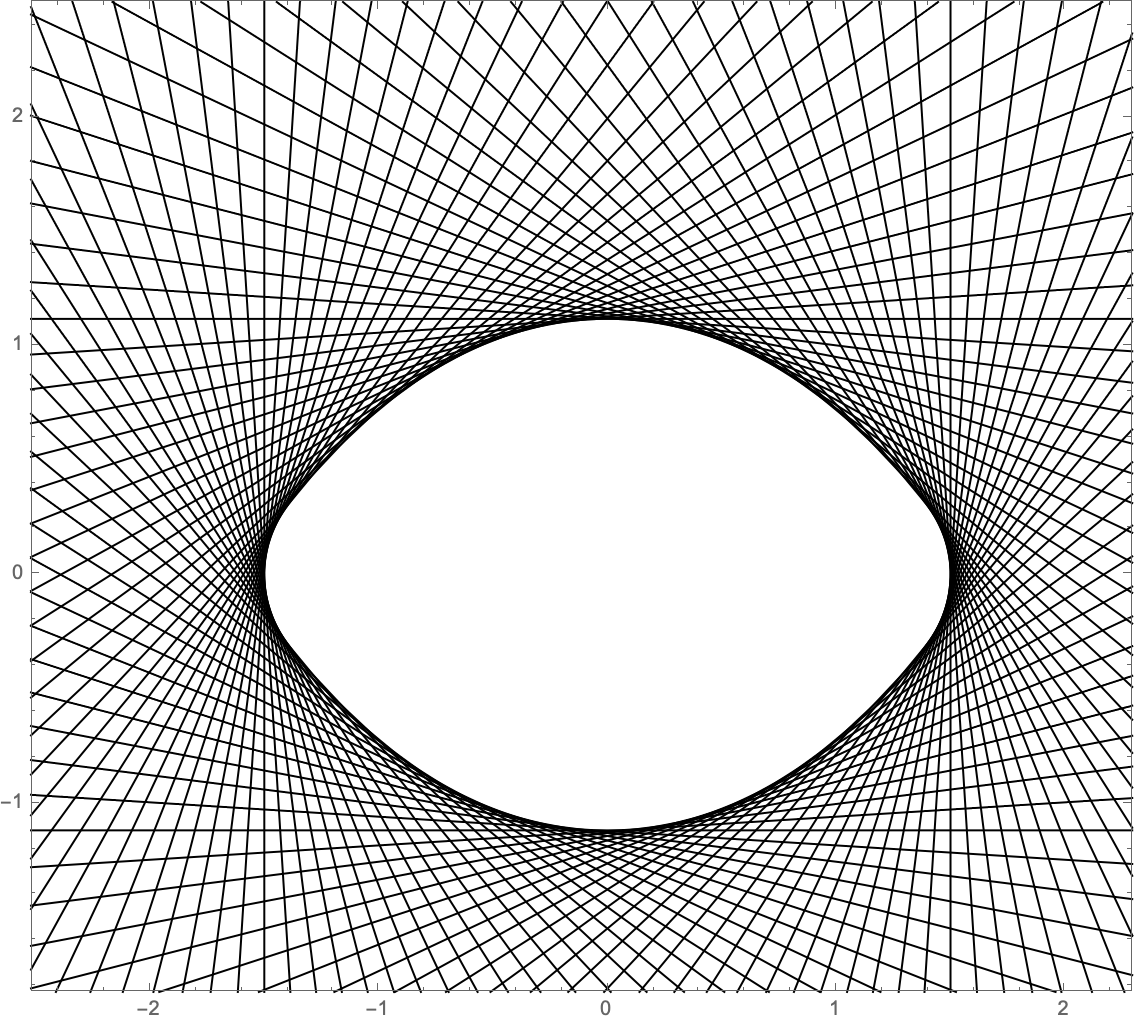}
    \caption{Supporting lines of $W(F_{I})$.}
    \label{fig1}
\end{figure}

\section{From supporting lines to the point description} \label{s:points} 

With formulas \eqref{lmax} describing the supporting lines, we can derive the point equation for the envelope curve explicitly. The respective standard procedure is to exclude $\theta$ from the system of equations 
\eq{ff} 	f(x,y,\theta)=0, \quad f'_\theta(x,y,\theta)=0, \en 
where the former represents the family of lines \eqref{sup} in a parametric form. In accordance with \eqref{lmax}, there are two cases to consider, depending on the relation between $\abs{\cos\theta}$ and the value of $(\sqrt{4+r^2}-r)/2$. 
\begin{prop} \label{th:smalltheta}Let $\abs{\cos\theta}>(\sqrt{4+r^2}-r)/2$. The respective portion of the envelope of \eqref{sup} 
consists of arcs of the circles of radius $r$ centered at $\pm 1$.
\end{prop} 
\begin{proof}Plugging $\lambda_{\max}(\theta)$ from the upper line of \eqref{lmax} into \eqref{sup} yields 
\[ f(x,y,\theta)=(x\mp 1)\cos\theta+y\sin\theta-r.\]  
From here, 
\[ f'_\theta =-(x\mp 1)\sin\theta+y\cos\theta, \]
and solving the respective system \eqref{ff} yields $x\mp 1=r\cos\theta, y=r\sin\theta$.  \end{proof} 

The remaining case is more involved. 
\begin{prop} \label{th:bigtheta} Let $\abs{\cos\theta}<(\sqrt{4+r^2}-r)/2$. The respective portion of the envelope of \eqref{sup} lies on the algebraic curve defined by \eqref{arc}. \end{prop} 
\begin{proof}With $\lambda_{\max}$ given by the second line in \eqref{lmax}, the system \eqref{ff} (at least, in the first quadrant --- which is sufficient for our purposes) takes the form 
\[ \begin{cases} x\cos\theta+y\sin\theta=\sqrt{(r/\sin\theta)^2+1}/2, \\ 
x\sin\theta-y\cos\theta=\frac{r^2\cos\theta}{2\sin^3\theta\sqrt{(r/\sin\theta)^2+1}}.\end{cases} \]
Replacing this system by the square of its first equation and its product with the second yields 
\[ \begin{cases} (x\cos\theta+y\sin\theta)^2=\left((r/\sin\theta)^2+1\right)/4, \\ 
	(x\sin\theta-y\cos\theta)(x\cos\theta+y\sin\theta)=r^2\cos\theta/\sin^3\theta.\end{cases} \]
Using the universal trigonometric substitution $t=\tan(\theta/2)$ turns the latter system into
 \eq{tp}
	\begin{cases}
		&-r^2 t^{10}-3 r^2 t^8-2 t^6 \left(r^2-8 x^2+8 y^2\right)+2 t^4 \left(r^2-8 x^2+8 y^2\right)\\& +3 r^2 t^2+r^2+8 t^7 x y-48 t^5 x y+8 t^3 x y = 0, 
		\\
		&-r^2 t^8-4 t^6 \left(r^2-x^2+1\right)-2 t^4 \left(3 r^2+4 x^2-8 y^2+4\right)\\& -4 t^2 \left(r^2-x^2+1\right)-r^2-16 t^5 x y+16 t^3 x y = 0.
	\end{cases}
\en
Equation \eqref{arc} is nothing but the resultant of \eqref{tp} considered as the system of two polynomials in $t$, obtained with the aid of \textit{Mathematica}.
\end{proof} 

\begin{figure}[H]
    \centering
    \includegraphics[width=0.8\textwidth]{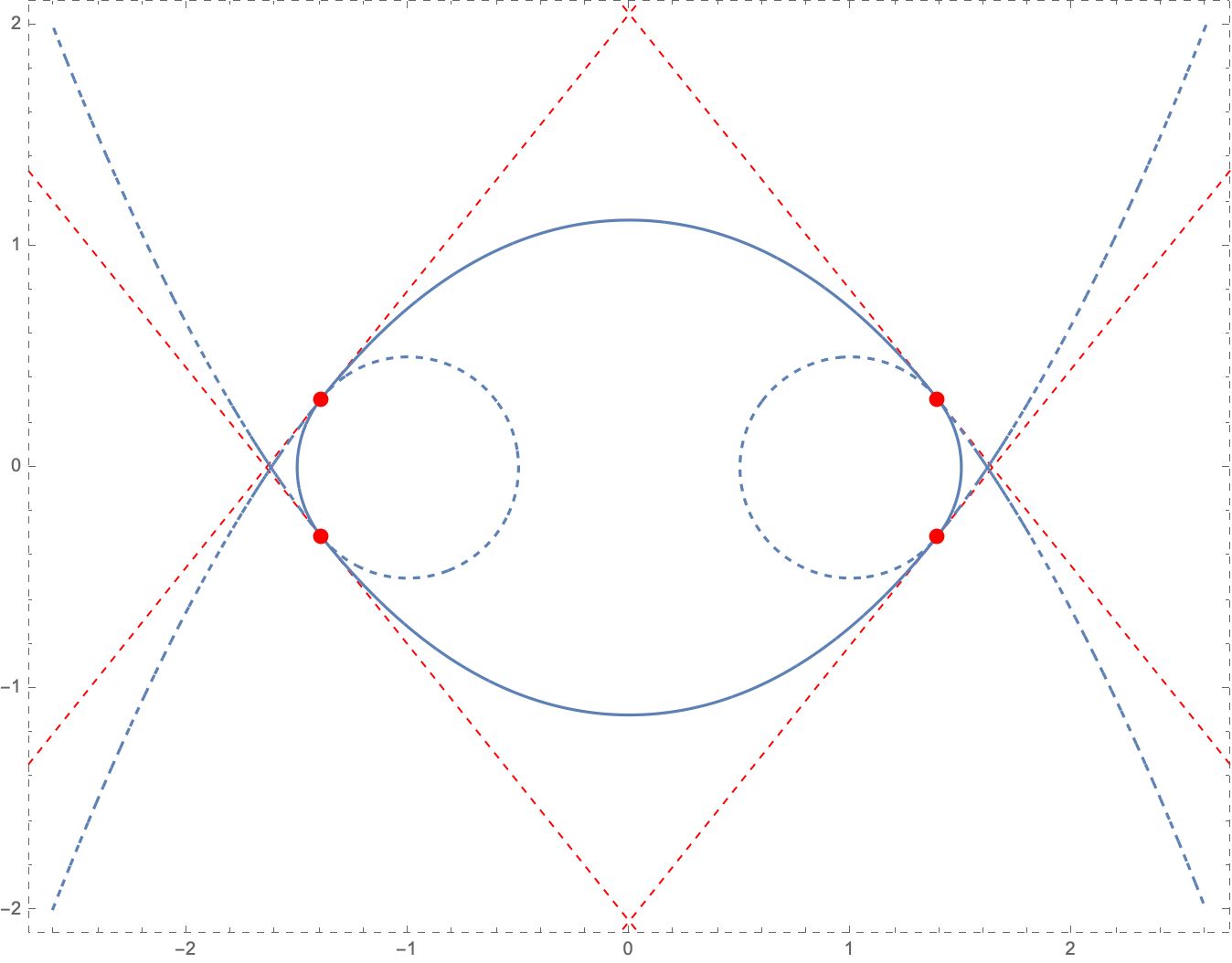}
    \caption{Boundary of $W(F_{I})$ (solid blue curve).}
    \label{fig2}
\end{figure}
In Fig.~\ref{fig2}, the circles and ``parabola-like'' curves containing the arcs of $\partial W(F_I)$ are plotted as blue dashed lines; red dashed lines are the supporting lines of $W(F_I)$ at the switching points (red dots) between the two types of boundary arcs.

\providecommand{\bysame}{\leavevmode\hbox to3em{\hrulefill}\thinspace}
\providecommand{\MR}{\relax\ifhmode\unskip\space\fi MR }
\providecommand{\MRhref}[2]{%
	\href{http://www.ams.org/mathscinet-getitem?mr=#1}{#2}
}
\providecommand{\href}[2]{#2}

\end{document}